\documentclass[11pt]{amsart}
\usepackage{graphicx}
\usepackage{amsfonts}
\usepackage{latexsym}
\usepackage{epsfig}

\newcommand{\beql}[1]{\begin{equation}\label{#1}}
\newcommand{\eeq}{\end{equation}}
\newcommand{\comment}[1]{}

\newcommand{\Set}[1]{{\left\{{#1}\right\}}}

\newcommand{\RR}{{\mathbb R}} 

\newcommand{\ZZ}{{\mathbb Z}}

\newcommand{\ft}[1]{\widehat{#1}}

\newcounter{open}
\newcounter{dfn}
\def\thedfn{\arabic{dfn}}

\newcounter{obs}
\def\theobs{\arabic{obs}}

\newcounter{thm}
\newcounter{othm}
\def\theothm{\Alph{othm}}

\newcounter{mysec}
\newcounter{mysubsec}[mysec]
\newtheorem{theorem}{Theorem}

\begin{document}

\title{Filling a box with translates of two bricks}

\author{Mihail N. Kolountzakis}

\thanks{
Supported in part by European Commission IHP Network HARP
(Harmonic Analysis and Related Problems),
Contract Number: HPRN-CT-2001-00273 - HARP.
}

\address{
School of Mathematics, Georgia Institute of Technology,
686 Cherry Street NW, Atlanta, GA 30332, United States, and\\
Department of Mathematics, University of Crete, Knossos Ave.,
714 09 Iraklio, Greece.
}

\email{kolount@member.ams.org}

\date{September 2004}

\begin{abstract}
We give a new proof of the following interesting fact recently proved by Bower and Michael
\cite{bower-michael}: if a $d$-dimensional rectangular box can be tiled using
translates of two types of rectangular bricks,
then it can also be tiled in the following way.
We can cut the box across one of its sides into two boxes, one of which can be tiled
with the first brick only and the other one with the second brick.
Our proof relies on the Fourier Transform.
We also show that no such result is true for three, or more, types of bricks.
\end{abstract}

\maketitle

Suppose we have at our disposal two types of $d$-dimensional rectangles (bricks),
type A with dimensions $(a_1,\ldots,a_d)$ and type B with dimensions
$(b_1,\ldots,b_d)$. 
We want to use translates of such bricks to fill completely, and with no overlaps,
a given $d$-dimensional rectangular box.
We then say that these two bricks tile our box by translations.

Bower and Michael \cite{bower-michael} recently showed the following
nice result.
A {\em hyperplane cut} is a seperation of an axis-aligned box in $d$ dimensions
using a hyperplane of the type $x_j=\alpha$,
for some $j=1,\ldots,d$ and some $\alpha \in \RR$.
A hyperplane cut separates such a box into two rectangular boxes
(all rectangles that appear in this note are axis-aligned).
\begin{theorem}\label{th:main}
{\rm (Bower and Michael \cite{bower-michael})}
If two bricks, of types A and B,  tile a box $Q$ (in dimension $d\ge 1$) by translations then
we can split $Q$ into two other boxes $Q_a$ and $Q_b$ using
a hyperplane cut, such that $Q_a$ can be tiled using translates of type A bricks
only and $Q_b$ can be tiled using translates of type B bricks only.
\end{theorem}
(For $d=1$ the result is obvious.)

The purpose of this note is to give a short proof of this fact using the Fourier
Transform, a very natural tool for this problem, as will become apparent.

\begin{figure}[htb]
\center{
\leavevmode
\epsfbox{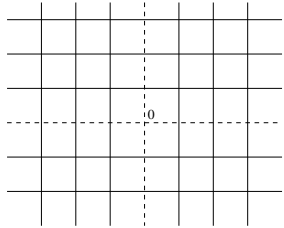}
\caption{The zeros (solid lines) of the Fourier Transform of a rectangle in 2 dimensions}
\label{fig:zeros}
}
\end{figure}

Indeed, suppose that $A = (-a_1/2,a_1/2)\times\cdots\times(-a_d/2,a_d/2)$ and
$B = (-b_1/2,b_1/2)\times\cdots\times(-b_d/2,b_d/2)$ are the two bricks and
$\Lambda_a, \Lambda_b$
are two finite subsets of $\RR^d$ which represent the translations of $A$ and $B$
that make up our box $Q = (-1/2,1/2)^d$
(as we may clearly assume without loss of generality).
In other words
\beql{tiling}
\sum_{\lambda \in \Lambda_a} \chi_A(x-\lambda) +
\sum_{\lambda \in \Lambda_b} \chi_B(x-\lambda) =
  \chi_Q(x),\ \ \mbox{a.e.\ $x\in\RR^d$}.
\eeq
The definition of the Fourier Transform $\ft{f}$ of a function $f \in L^1(\RR^d)$
that we use is
$$
\ft{f}(\xi) = \int_{\RR^d} f(x) \exp(-2\pi i \xi\cdot x)\,dx.
$$
Taking the Fourier Transform of both sides of \eqref{tiling}
we get
\beql{tiling-ft}
\phi_a(\xi)\ft{\chi_A}(\xi) + \phi_b(\xi)\ft{\chi_B}(\xi) = \ft{\chi_Q}(\xi),
\eeq
where $\phi_a(\xi) = \sum_{\lambda\in\Lambda_a} \exp(2\pi i \lambda\cdot x)$,
$\phi_b(\xi) = \sum_{\lambda\in\Lambda_b} \exp(2\pi i \lambda\cdot x)$, are 
trigonometric polynomials.
Simple calculation shows that the Fourier Transform of the indicator
function of the box $C = (-c_1/2,c_1/2)\times\cdots\times(-c_d/2,c_d/2)$ is
\beql{box-ft}
\ft{\chi_C}(\xi) = \prod_{j=1}^d \frac{\sin(c_j\xi_j)}{\xi_j},
\eeq
whose zero set $Z(\ft{\chi_C})$ consists of all points $\xi$ with at least one coordinate
$\xi_j$ being a non-zero multiple of $c_j^{-1}$.
This set may be viewed as a collection of $d$ sets of hyperplanes, with the hyperplanes
in the $j$-th set being parallel to the hyperplane $\xi_j=0$ and spaced at regular
intervals $c_j^{-1}$, with the exception of the hyperplane $\xi_j=0$
itself (see Figure \ref{fig:zeros}).
\epsfysize=3cm

Therefore the zero set of the right hand side of \eqref{tiling-ft}
is the set
\beql{rhs-zeros}
Z = Z(\ft{\chi_Q}) = \Set{\xi\in\RR^d: \xi_j \in \ZZ\setminus\Set{0},\mbox{\ for some $j=1,\ldots,d$}}.
\eeq
The key observation is the following: for any choice of different $i$ and $j$ from
the numbers $1,\ldots,d$ at least one of $a_i^{-1}$ and $b_j^{-1}$ is an integer.
For, assuming otherwise, the hyperplanes $\xi_i=a_i^{-1}$ and $\xi_j = b_j^{-1}$ would
be part of the zeros sets of the first and second term in the left hand side of \eqref{tiling-ft}
respectively.
But the intersection of these hyperplanes, on which set the left hand side vanishes,
contains points not in the set $Z$ of \eqref{rhs-zeros}, a contradiction.

Finally, if the numbers $a_1^{-1},\ldots,a_d^{-1}$ are all integers then brick $A$ can tile
$Q$ alone and there is nothing to prove. So we may assume that one of them is not an integer,
say $a_1^{-1} \notin \ZZ$.
By choosing $i=1$ and $j=2,3,\ldots,d$ in turn,
and using our key observation above, we deduce that all $b_j^{-1}$, $j=2,3,\ldots,d$, are integers.
For the same reason as before we can also assume that $b_1^{-1}$ is not an integer (otherwise
brick $B$ can tile alone), which in turn shows that all $a_j^{-1}$, $j=2,3,\ldots,d$,
are integers.
Hence the face of each brick parallel to the $x_1=0$ hyperplane can tile the corresponding face of $Q$.

On the other hand, by the assumed tiling of $Q$ by translates of bricks $A$ and $B$ it follows,
by looking along the first coordinate axis, that $1 = m a_1 + n b_1$ for some nonnegative integers
$m$ and $n$. Split then the box $Q$ by the hyperplane $x_1 = -1/2 + m a_1$ into two boxes
of dimensions $m a_1\times 1 \times \cdots \times 1$ and $n b_1\times 1 \times \cdots \times 1$.
The first box can tiled by brick $A$ by simply tiling its $1\times\cdots\times1$ face and repeating
this $m$ times. The second box can be tiled similarly by box $B$, as we had to show.

\noindent
{\bf An example.}
Let us observe that there is no generalization of this result to three or more bricks.
That is, there are boxes which admit tilings with translates of three types of bricks, but which
cannot be split into two parts using a hyperplane cut so that each of these parts
can be tiled with a proper subset of the available types of bricks.
It is enough to give an example in two dimensions, as any such example can be transformed
to one in dimension $d>2$ be considering all bricks to have their last $d-2$ coordinates equal to
1, and considering the $d$-dimensional tiling that arises by one layer of the two-dimensional
example.
\epsfysize=3cm
\begin{figure}[htb]
\center{
\leavevmode
\epsfbox{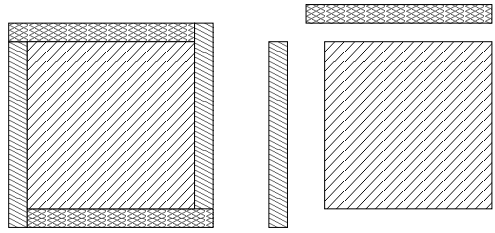}
\caption{A tiling of a rectangle (left) with three types of bricks (right).}
\label{fig:three}
}
\end{figure}

To see a two-dimensional example take $R$ much larger than $1$ and use the three brick types
$1\times R$, $R \times 1$ and $(R-1)\times(R-1)$. With these we can tile a $(R+1)\times(R+1)$ box as
shown in Figure \ref{fig:three}.
But the box cannot be split into two boxes using a hyperplane cut, each of which can
be tiled using a proper subset of the available brick types.
This can be verified by examining the few possibilities.


\end{document}